\documentclass[12pt,a4paper]{amsart}
\usepackage{amsfonts,amssymb,amscd,amstext,mathrsfs}
\theoremstyle{plain}

\usepackage[colorlinks]{hyperref}
\usepackage{enumerate, amssymb}

\advance\hoffset-5mm \advance\textwidth40mm
\input{diagrams.tex}

\def\bdi{\begin{diagram}}
\def\edi{\end{diagram}}

\theoremstyle{plain}

\newtheorem{thm}{Theorem}[section]
\newtheorem{cor}[thm]{Corollary}
\newtheorem{lem}[thm]{Lemma}
\newtheorem{prop}[thm]{Proposition}
\theoremstyle{definition}
\newtheorem{defi}[thm]{Definition}
\newtheorem{defis}[thm]{Definitions}
\newtheorem{conj}[thm]{Conjecture}
\newtheorem{conv}[thm]{Convention}
\newtheorem{nota}[thm]{Notation}
\newtheorem{rem}[thm]{Remark}
\newtheorem{rems}[thm]{Remarks}
\newtheorem{exa}[thm]{Example}
\newtheorem{exas}[thm]{Examples}
\newtheorem{prob}[thm]{Problem}
\newtheorem{probs}[thm]{Problems}
\newtheorem{ques}[thm]{Question}
\newtheorem{sit}[thm]{}

\newcommand{\Ker}{ \operatorname{{\rm Ker}}}

\newcommand{\Aut}{ \operatorname{{\rm Aut}}}

\newcommand{\SL}{ \operatorname{{\bf SL}}}

\def\deg{\mathop{\rm deg}}

\renewcommand{\epsilon}{\varepsilon}

\def\and{\quad\mbox{and}\quad}

\newcommand{\C}{\ensuremath{\mathbb{C}}}

\newcommand{\N}{\ensuremath{\mathbb{N}}}

\newcommand{\bX}{{\bar X}}

\newcommand{\tX}{{\tilde X}}
\newcommand{\tY}{{\tilde Y}}
\newcommand{\tZ}{{\tilde Z}}

\newcommand{\bZ}{{\bar Z}}

\newcommand{\cF}{{\ensuremath{\mathcal{F}}}}

\newcommand{\cA}{{\ensuremath{\mathcal{A}}}}

\newcommand{\cU}{{\ensuremath{\mathcal{U}}}}

\newcommand{\p}{\partial}
\newcommand{\de}{\delta}
\newcommand{\id}{{\rm id}}

\renewcommand{\rho}{\varrho}

\def\bals#1\eals{\begin{align*}#1\end{align*}}
\def\bal#1\eal{\begin{align}#1\end{align}}

\def\SAut{\mathop{\rm SAut}}

\def\PP{{\mathbb P}}

\renewcommand{\phi}{\varphi}

\newcommand{\bnum}{\begin{enumerate}}
\newcommand{\enum}{\end{enumerate}}
\renewcommand{\emptyset}{\varnothing}

\newcommand{\brem}{\begin{rem}}
\newcommand{\brems}{\begin{rems}}
\newcommand{\erem}{\end{rem}}
\newcommand{\erems}{\end{rems}}
\newcommand{\bprob}{\begin{prob}}
\newcommand{\eprob}{\end{prob}}
\newcommand{\bprobs}{\begin{probs}}
\newcommand{\eprobs}{\end{probs}}
\newcommand{\bques}{\begin{ques}}
\newcommand{\eques}{\end{ques}}
\newcommand{\bexa}{\begin{exa}}
\newcommand{\bexas}{\begin{exas}}
\newcommand{\eexa}{\end{exa}}
\newcommand{\eexas}{\end{exas}}
\newcommand{\bdefi}{\begin{defi}}
\newcommand{\edefi}{\end{defi}}
\newcommand{\bdefis}{\begin{defis}}
\newcommand{\edefis}{\end{defis}}
\newcommand{\bcor}{\begin{cor}}
\newcommand{\ecor}{\end{cor}}
\newcommand{\blem}{\begin{lem}}
\newcommand{\elem}{\end{lem}}
\newcommand{\bconv}{\begin{conv}}
\newcommand{\econv}{\end{conv}}
\newcommand{\bconj}{\begin{conj}}
\newcommand{\econj}{\end{conj}}
\newcommand{\bprop}{\begin{prop}}
\newcommand{\eprop}{\end{prop}}
\newcommand{\bthm}{\begin{thm}}
\newcommand{\ethm}{\end{thm}}
\newcommand{\bnota}{\begin{nota}}
\newcommand{\enota}{\end{nota}}
\newcommand{\bsit}{\begin{sit}}
\newcommand{\esit}{\end{sit}}
\newcommand{\be}{\begin{equation}}
\newcommand{\ee}{\end{equation}}
\newcommand{\bproof}{\begin{proof}}
\newcommand{\eproof}{\end{proof}}
\def\ba{\begin{array}}
\def\ea{\end{array}}

\thanks{The
second author was  partially supported by Schweizerische
Nationalfonds grants No. 200020-134876/1 and 200021-140235/1 and the third author was supported by Australian Research Council grants DP120104110 and DP150103442.}

\begin{document}
\title[On subelliptic manifolds]{On subelliptic manifolds}

\author{Shulim Kaliman, 
Frank Kutzschebauch and Tuyen Trung Truong}
\address{
University of Miami, Coral Gables, FL 33124, USA}
\email{kaliman@math.miami.edu}
\address{Department of Mathematics,
University of Bern, Bern, Switzerland}
\email{frank.kutzschebauch@math.unibe.ch	}
\address{School of Mathematical Sciences, University of Adelaide, Adelaide SA 5005, Australia}
\email{tuyen.truong@adelaide.edu.au}

\date{\today}
\maketitle

\hspace{6.5cm}  {\em Dedicated to Mikhail Zaidenberg} 

\hspace{6.5cm} {\em on the occasion of his 70-th birthday}




\begin{abstract} 
A smooth complex quasi-affine algebraic variety $Y$ is flexible if its special group $\SAut (Y)$ of automorphisms (generated by the elements of one-dimensional unipotent subgroups of $\Aut (Y)$) acts transitively on $Y$. An irreducible algebraic manifold $X$ is locally stably flexible if it is the union $\bigcup X_i$ of a finite number of Zariski open sets, each $X_i$ being quasi-affine, so that there is a positive integer $N$ for which $X_i\times \mathbb{C}^N$ is flexible for every $i$. 
The main result of this paper is  that the blowup of a locally stably flexible manifold at a smooth algebraic submanifold (not necessarily equi-dimensional or connected) is subelliptic, and hence Oka. This result is proven as a corollary of some general results concerning the so-called $k$-flexible manifolds. 
 
\end{abstract}

\thanks{
{\renewcommand{\thefootnote}{} \footnotetext{ 2010
\textit{Mathematics Subject Classification:}
14R20,\,32M17.\mbox{\hspace{11pt}}}}

{\footnotesize \tableofcontents}

\section*{Introduction} The notion of a subelliptic manifold (i.e. a manifolds which admits a dominating family of  sprays) was introduced by Forstneri\v c in \cite{For1},
inspired by hints from Gromov in \cite{Gr2}. It is a natural generalization to the stronger condition of admitting a single dominating spray, called elliptic. The importance of the notion 
of subellipticity is that as in the case of elliptic manifolds it implies all Oka properties.  In other words such a subelliptic manifold is $X$ is an Oka manifold as proven by Forstneri\v c in  \cite{For1}. In particular being an Oka manifold implies that every holomorphic map from a convex domain $K$
in $\C^n$ into $X$ can be approximated (in the compact-open topology) by a holomorphic map from $\C^n$ to $X$. 
Needless to say that this leads to many remarkable consequences (e.g., see \cite{For}). 
On the other hand, having the same consequences, subellipticity is easier to establish than ellipticity,  which  is exemplified by the main results of the present paper.

The simplest example of an elliptic manifold is, of course,  the Euclidean space $\C^n$ itself. 
Furthermore, Gromov proved ellipticity in the case of the complement to a subvariety
of codimension at least 2 in $\C^n$. Any algebraic manifold which is locally isomorphic to such complements (resp. $\C^n$)
is called a manifold of class $\cA$ (resp. $\cA_0$) (\cite[Definition 6.4.5. ]{For}, \cite[Remark 3] {LT}) Since  in the algebraic case subellipticity turns out to be a local property we see that
a manifold of class $\cA$ is always subelliptic. Gromov observed also the following.

\bprop\label{in.p1} Let $X$ be a complex manifold of class $\cA_0$ and $Y$ be the result of blowing $X$ up at a finite number of
points. Then $Y$ is also a manifold of class $\cA_0$ and, therefore, subelliptic.
\eprop

For example this yields subellipticily of compact rational surfaces (see \cite[Corollary 6.4.8 ]{For}).  

%

There were no analogs of Proposition \ref{in.p1} until the recent paper of L\'arusson and the third author \cite{LT} who proved the following.

\bthm\label{in.t2} Let $X$ be an algebraic manifold of class $\cA$ and $\pi : \tX \to X$ be the blowing up of $X$ along a smooth algebraic
(not necessarily connected) submanifold of codimension at least 2. Then $\tX$ is subelliptic.

\ethm

The proofs in \cite{LT} made use crucially of the fact that $\mathbb{C}^n$ has a lot of automorphisms, as first pointed out by Gromov and generalised by Winkelmann. This last property is shared by the so-called flexible manifolds, extensively-studied in affine algebraic geometry. Recall that one of equivalent definitions states that  a smooth complex quasi-affine algebraic variety $X$ of dimension at least 2 is flexible if its special group $\SAut (X)$ of automorphisms (generated by the elements of one-dimensional unipotent subgroups of $\Aut (X)$) acts transitively on $X$.
It is easy to establish that flexible manifolds are algebraically subelliptic (and even algebraically elliptic). Furthermore, there is no need to discuss complements to subvarieties
of codimension at least 2 in flexible manifolds because such complements are again flexible (\cite{FKZ-GW}).
This observation was a strong indication for us that the above construction can survive replacement of Euclidean spaces by flexible manifolds. This is in fact true, and we can actually prove the same result for a more general class of manifolds which we are going to define next. 

{\bf Definition.} An irreducible algebraic manifold $X$ is locally stably flexible, if it is the union $\bigcup X_i$ of a finite number of Zariski open sets, each $X_i$ being quasi-affine, so that there is a positive integer $N$ for which $X_i\times \mathbb{C}^N$ is flexible for every $i$. 

Our main result is the following.

\bthm\label{i.t1} Let $X$ be a locally stably flexible manifold. Suppose that $\pi : \tX \to X$ is the blowing up of $X$ along a smooth algebraic submanifold $Z$, not necessarily equi-dimensional  or connected. Then $\tX$ is algebraically subelliptic. 
\ethm

After the Theorems \ref{in.t2} and \ref{i.t1}, it is natural to ask the following questions concerning the behavior under blowups of various classes of algebraic manifolds of interest in affine geometry and Oka theory. 

{\bf Q1.} Is the class of algebraically subelliptic manifolds preserved by blowups? We note that by Gromov's results the class $\mathcal{A}$ (which is a smaller class of manifolds) and by results in \cite{LT} the class of strongly algebraically dominated manifolds (which is a bigger class of manifolds) are both preserved by blowups.  

{\bf Q2.} Is the class of locally stably algebraic manifolds preserved by blowups?



The remaining of this paper is organized as follows.
In Section 1 we remind the notions of sprays and subellipticity and establish some simple facts which are immediate consequences of the results
presented in \cite{For}. In Section 2 we describe technique developed for flexible manifolds in \cite{FKZ-GW} and \cite{AFKKZ} and prove a non-trivial
fact (Theorem \ref{fm.t1}) heavily based on \cite{FKZ-GW}. This theorem deals with a partial quotient morphism $\rho : X \to Q$ of a flexible manifold
$X$ with respect to some $G_a$-action. It establishes that up to an automorphism of $X$ for every closed submanifold $Z$ of $X$ (of codimension at least 2)
and every $z \in Z$ we can suppose that $\rho|_Z : Z \to \rho (Z)$ is a local isomorphism over a Zariski neighborhood of $\rho (z) \in \rho (Z)$. 
In Section 3 we prove simple facts which, in particular, include ellipticity of flexible manifolds.
Section 4 is devoted to technique of affine modifications which can be mostly found in \cite{KZ}. It is necessary because in
the above notations $X$ is an affine modification of $Q \times \C$. It turns out that we need to present $X$ as a more refined
affine modification for which we introduce in Section 5 the notion of $k$-flexibility. Namely, $X$ is $k$-flexible if for some dominant morphism $\tau : X \to P$ there is a Zariski dense open subset
$P_0$ of $P$ for which $\tau^{-1}(P_0)$ is 
 isomorphic to $P_0\times \C^k$ which implies that $X$ is an affine modification
of $P\times \C^k$. \footnote{It is worth mentioning that most of flexible manifolds are $k$-flexible for some $k\geq 2$. Actually, starting from dimension 3
we do not know examples of flexible manifolds that are not at least 2-flexible.}
With all preparations done we obtain our main theorems in Section 6. 

{\bf Acknowledgments.} The authors would like to thank Finnur L\'arusson for his helpful comment concerning the descent property of algebraic subellipticity. 
 Most of this work was done during a stay of the second author at the University of Miami
and he thanks this institution for its hospitality and excellent working conditions. It was also partially done while the second and third authors were attending the program "Workshop on higher algebraic geometry, holomorphic dynamics and their interaction" at the Institute for Mathematical Sciences, National University of Singapore in January 2017, and we thank the institution and the organizers for their hospitality and financial support. 

\section{Sprays and subellipticity}

Let us remind some definitions which can be found in \cite{For}.

\bdefi\label{sas.d1} (i) A holomorphic vector bundle $p : E \to X$ over a complex manifold $X$ is called a spray is there exists
a holomorphic map $s : E \to X$ such that for every point $y$ in the zero section $S$ of $E$ one has $s (y)=p(y):=x$.
That is, a spray is a triple $(E, p, s)$.

(ii) A spray is called dominating if for every $y \in S$ one has ${\rm d}s (T_yp^{-1}(x))=T_xX$.

(iii) A family of sprays $\{E_i,p_i,s_i\}_{i=1}^m$ on $X$ is called dominating if for every $y \in S$
$${\rm d}s_1 (T_yp_1^{-1}(x)+ {\rm d}s_2 (T_yp_2^{-1}(x)+\ldots + {\rm d}s_m (T_yp_m^{-1}(x)=T_xX.$$

(iv) A complex manifold $X$ is called elliptic (resp. subelliptic) if it admits a dominating holomorphic spray (resp. a dominating family
of holomorphic sprays).

(v) We say that a spay  $(E, p, s)$ is of rank $k$ if the rank of the vector bundle $p : E \to X$ is $k$
\edefi

\bconv\label{sas.con1} From now on we consider only {\bf algebraic} sprays $(E,p,s)$ on algebraic complex manifolds which means that the vector
bundle $p : E \to X$ is algebraic and the map $s: E \to X$ is algebraic. We omit this adjective ``algebraic" below.

\econv

Under this convention the following definition makes sense.

\bdefi\label{sas.d2}
(a) Let $X_0$ be a nonempty Zariski open subset of a complex algebraic manifold $X$. An algebraic vector bundle $p : E \to X_0$ is called a spray on $X_0$ with values in $X$ if  there exists
a holomorphic map $s : E \to X$ such that for every point $x\in X_0$ in the zero section $S_0$ of $E$ one has $s (x)=x$.

(b) Let $s': E\to X_0$ be another spray on $X_0$ with values in $X$ where $p: E \to X$ is the same vector bundle as in (a). We say that it is equivalent to the spray $s$ from (a) if for general points
$y \in S_0$ and  $x=p(y)$ there is a linear automorphism $\lambda$ of the fiber $E_x=p^{-1}(x)$ for which 
\be\label{sas.e1} s\circ \lambda|_{E_x}= s'|_{E_x} \ee 

(c) The notion of a dominating spray on $X_0$ with values in $X$ is described exactly as in Definition \ref{sas.d1} with $y$ running over
a section $S_0$ of $p : E \to X_0$. In the same fashion we deal with a dominating family of sprays on $X_0$ with values in $X$.

\edefi

Convention \ref{sas.con1} enables us to use the following facts (see, the proof of \cite[ Theorem 6.4.2]{For}) again and again.

\bprop\label{sas.p1} Let $s : E \to X$ be a spray on $X_0$ with values in $X$ as in Definition \ref{sas.d2}. 
Then there exists an equivalent spray $s': E\to X$ on $X_0$ with values in $X$ such that it extends to a spray on $X$.
Furthermore, if $X\setminus X_0$ is
a principal divisor  then this spray $s'$ can be chosen so that equation (\ref{sas.e1}) holds  for every $x \in X_0$.
\eprop

\bcor\label{sas.c1} Let $\{U_i \}_i$ be a cover of a complex algebraic manifold $X$ by Zariski open sets such that for every $i$
there is a dominating family of sprays on $U_i$ with values in $X$. Then there is a dominating family of sprays on $X$.

\ecor

\bcor\label{sas.c2} Let $s : E \to X$ be a spray on $X$ as in Definition \ref{sas.d1} and let
$\varphi : X \to Y$ be a birational morphism which yields an isomorphism between Zariski open subsets $X_0 \subset X$ and $Y_0 \subset Y$,
i.e. one has the following commutative diagram
\[  \begin{array}{ccc} E|_{X_0} \, \, \, &  \stackrel{{\psi}}{\rightarrow} &F \\
\, \, \, \, \downarrow {p|_{X_0}}  & & \, \, \, \, \downarrow { q} \\
X_0\, \, \,\, \, \,&  \stackrel{{\varphi|_{X_0}}}{\rightarrow} & Y_0.
\end{array} \]
of isomorphic vector bundles. Then $r=\varphi \circ s \circ \psi^{-1}: F \to Y$ is a spray on $Y_0$ with values in $Y$.
Furthermore, if the complement of $Y_0$ in $Y$ is a principal divisor there  is  an equivalent spray $r': F \to Y$ extendable to a spray on $Y$
and such that $r(q^{-1}(y))=r'(q^{-1}(y))$ for every $y\in Y_0$.  
\ecor

\bcor\label{sas.c4} The class of manifolds admitting sprays is closed with respect to the  procedure of blowing down.
\ecor

\section{Flexible manifolds}

Recall the following facts which can be found in \cite{AFKKZ}.

\bdefi\label{fm.d1}  (1) A derivation $\sigma$ on the ring $A$ of regular functions on a quasi-affine algebraic manifold $X$ is called locally nilpotent
if for every $0\ne a \in A$ there exists a natural $n$ for which $\sigma^n (a)=0$. For the smallest $n$ with this property one defines
the degree of $a$ with respect to $\sigma$ as $\deg_\sigma a=n-1$. This derivation can be viewed as a vector field on $X$ which
we also call locally nilpotent. The phase flow of this vector field is an algebraic $G_a$-action on $X$, i.e. the action of the group $\C_+$
of complex numbers with respect to addition 
which can be viewed as a one-parameter unipotent group $U$ in the group $\Aut (X)$ of all algebraic automorphisms of $X$.
In fact, every $G_a$-action is generated by a locally nilpotent vector field (e.g, see \cite{Fre}).

(2) A quasi-affine manifold $X$ is called flexible if for every $x \in X$ the tangent space $T_xX$ is spanned by the tangent vectors
to the orbits of one-parameter unipotent subgroups of $\Aut (X)$ through $x$.

(3) The subgroup $\SAut (X)$ of $\Aut X$ generated by all one-parameter unipotent subgroups is called special.
\edefi

We have the following \cite{AFKKZ},  \cite{FKZ-GW}. 

\bprop\label{fm.p1} For every irreducible  quasi-affine  algebraic variety $X$ the following are equivalent

{\rm (i)} the special subgroup $\SAut (X)$ acts transitively on $X_{reg}$;

{\rm (ii)}  the  special subgroup $\SAut (X)$ acts infinitely transitively on $X_{reg}$ (i.e. for every natural $m$
the action is $m$-transitive);

{\rm (iii)} $X_{reg}$ is flexible.

\eprop

By the Rosenlicht Theorem (see \cite[Theorem
2.3]{PV}) for $X$, $A$, and $U$ as in Definition \ref{fm.d1}
one can find a finite set of $U$-invariant functions $a_1,\ldots,a_m\in A$,
which separate general $U$-orbits in $X$. They generate a morphism $\rho : X \to Q$ into an affine
algebraic variety $Q$. Note that this set of invariant functions can be chosen so that $Q$ is normal (since $X$ is normal).

\bdefi  Such a morphism $\rho : X \to Q$ into a normal $Q$ will be called a partial quotient.
In the case when $a_1, \ldots , a_m$ generate the subring $A^U$ of $U$ invariant elements of $A$
such a morphism is called the categorical quotient.\footnote{However, in general $A^U$ is not finitely generated by the 
Nagata's example. That is, why, following \cite{FKZ-GW} we prefer to work with partial quotients.}
\edefi

The main aim of this section is the next theorem.

\bthm\label{fm.t1} Let $Z$ be a submanifold of codimension at least 2 in a flexible affine algebraic manifold $X$, and $\sigma$ be a nontrivial locally nilpotent
vector field on $X$.  Suppose that $\rho : X \to Q$ is a partial quotient morphism of the $G_a$-action associated with $\sigma$ such that $Q$ is a normal
variety. Then for every finite set $z_1, \ldots , z_m \in Z$ of distinct points one can find an automorphism $\alpha$ of $X$ 
such that for every $i=1, \ldots , m$ and
$\rho_\alpha = \rho \circ \alpha$ one has

{\rm (i)} the point $\rho_\alpha (z_i)$ is general (and, therefore, smooth) in $Q$,

{\rm (ii)} $\rho_\alpha (z_i)$  is contained in a smooth part of $\rho_\alpha (Z)$;

{\rm (iii)} the morphism $\rho_\alpha|_Z : Z \to \rho_\alpha (Z)$ is local embedding at $z_i$ (and, in particular,  $\rho_\alpha|_Z $ is birational).

Furthermore, let $\tau: Q\to P$ be a morphism such that $\dim P=\dim Z$. Then $\alpha$ can be chosen so that

{\rm (iv)}  $\tau|_{\rho_\alpha (Z)} : \rho_\alpha  (Z) \to P$
is \'etale at $\rho_\alpha (z_i)$ for every $i$.
\ethm

The proof of this fact is heavily based on the technique from \cite{FKZ-GW} and it requires some preparations, but
first let us extract some corollary using following notion introduced by Ramanujam \cite{Ra1}.

\bdefi\label{fm.d2}  Given irreducible algebraic varieties $X$ and $A$ and
a map $\varphi:A\to\Aut(X)$ we say that $(A,\phi)$
is an {\em algebraic family of automorphisms on $X$} if the induced map
$A\times X\to X$, $(\alpha,x)\mapsto \varphi(\alpha).x$, is a morphism.
\edefi

\brem\label{fm.r1}  Note that properties (i) -(iv) from Theorem \ref{} survive under a small perturbation of the automorphism $\alpha$
in an irreducible algebraic family of automorphisms $A$. This implies that they are valid for a Zariski open subset of $A$
(because of the algebraicity) and we have the following.
\erem

\bcor\label{fm.c1} Let Theorem \ref{fm.t1} hold for an automorphism $\alpha$ which is contained in an irreducible algebraic family of automorphisms $A$.
Then the statement remains valid if one replaces $\alpha$ by a general element of $A$.
\ecor

The first step in the proof of Theorem \ref{fm.t1} is the following.

\blem\label{fm.l1} Let the assumption of Theorem \ref{fm.t1} hold.
Then for every finite set $z_1, \ldots , z_m \in Z$ of distinct points one can find an automorphism $\alpha$ of $X$ 
such that for every $i=1, \ldots , m$ and
$\rho_\alpha = \rho \circ \alpha$ properties (i) and (iv) from Theorem \ref{fm.t1} are true. Furthermore,

{\rm (ii$'$)} $\rho |_{Z_\alpha} : Z_\alpha \to \rho (Z_\alpha)$ is a local embedding at every point $\alpha (z_i)$.

\elem

\bproof Choose general points $q_1, \ldots , q_m$ in $Q$ and general points $x_1, \ldots , x_m$ in $X$ for which $\rho (x_i)=q_i$.
Let $(v_{i,1}, \ldots , v_{i,n})$ (resp. $(u_{i,1}, \ldots , u_{i,n+1})$) be a local analytic coordinate system at $q_i \in Q$ (resp. $x_i \in X$) such that
$\rho^* (v_{i,j})=u_{i,j}$. By \cite[Theorem 4.14 and Remark 4.16]{AFKKZ} one can choose an automorphism $\alpha$ of $X$ such that
$\alpha (z_i)=x_i, \, i=1, \ldots, m$ and, furthermore, $\alpha(Z)$ is tangent to the subvariety $u_{i,n-k+1}=\ldots = u_{i, n+1}=0$ where $k+1$ is
the codimension of $Z$ in $X$. By construction $\alpha$ satisfies properties (i) and (ii$'$). 

For (iv) it suffices to require  that $(v_{i,1}, \ldots , v_{i,n-k})|_{\rho_\alpha (Z)}$ is a lift of a local analytic coordinate system on $P$ under $\tau$ and we are done.

\eproof

\bdefi\label{fm.d3} For every locally nilpotent vector fields $\sigma$ and each function $f \in \Ker \sigma$ from its kernel the field
$f\sigma$ is called a replica of $\sigma$. Recall that such a replica is automatically locally nilpotent.

\edefi

\bprop\label{fm.p1}{\rm (cf. \cite{FKZ-GW})} Let $\delta_0$ be a locally nilpotent vector
field on  a quasi-affine algebraic manifold $X$,
$\rho_0: X \to Q_0$ be an associated partial quotient morphism, $x$ be a general point of $X$, and $O_1$ be the orbit of $x$ under the phase flow of $\delta_0$.
Then  there exists a locally nilpotent vector field $\delta_1$ such that

{\rm ($\#$)}  for general points $x_1, \ldots , x_{n-1}\in O_1$ and the vectors $\delta_{1,x_1}, \ldots , \delta_{1,x_{n-1}}$ (which are the values of $\delta_1$ 
at these points) the vectors $v_1, \ldots , v_{n-1}$ form a basis of $T_{q_0}Q_0$ where $q_0=\rho_0(x)$ and $v_i ={\rm d}\rho_0 (\delta_{1,x_i})$.

Furthermore, let condition ($\#$) hold and $H$ be  the group of 
algebraic automorphisms of $X$ generated by the elements from the phase flows of $\delta_0$, $\delta_1$, and their replicas. 
Then the orbit of $x$ under the action of $H$  is Zariski open in $X$. 
\eprop

\bproof 
By \cite[Theorem 4.14 and Remark 4.16]{AFKKZ}
there exists an automorphism $\alpha \in \Aut (X)$ such that it fixes points $x_1, \ldots , x_{n-1}$ and for every $i$ the linear map ${\rm d}\alpha |_{T_{x_i}X}$ coincides
with a prescribed element of $\SL_{n-1}$. Furthermore, for every fixed $k \in \N$ (however large it is) we can require that the $k$-jet $\alpha_{x_i}^k$ of
$\alpha$ at $x_i$ coincides with such a linear part.
Hence choosing any locally nilpotent derivation $\delta_1$ for which every $\delta_{1,x_i} \ne 0$ we can achieve ($\#$) replacing $\delta_1$ by $\alpha_*(\delta_1)$.

For every nonzero locally nilpotent $\delta_0$ the statement  of \cite[Proposition 1.14]{FKZ-GW} yields the existence of $\delta_1$ for which the orbit of $x$ under the action of the group $H$ is open. 
However, the analysis of the proof of this
fact shows that any $\delta_1$ satisfying ($\#$) fits this purpose. 

\eproof

\brem\label{fm.r1}  (a) In fact, we can replace condition  ($\#$)  in Proposition \ref{fm.p1} with the  following.
For every point $y \in O_1$ denote by $y(t)$ the image of $y$  under the action of the phase flow of $\delta_0$ at time $t$.
Then we can require that

 ($\#'$)  for a general moment of time $t$ and general $y\in O_1$ the vector ${\rm d}\rho_0 (\delta_{1,x(t)})-{\rm d}\rho_0 (\delta_{1,y(t)})$
is general in $T_{q_0}Q_0$. 

The fact that  ($\#'$)  implies  ($\#$)  is clear. Indeed, otherwise
there is a proper subspace  $V \subset T_qQ_0$ such that for every point $w \in O_1$ one has
${\rm d}\rho_0 (\delta_{1,w}) \in V$. Then ${\rm d}\rho_0 (\delta_{1,x(t)})-{\rm d}\rho_0 (\delta_{1,y(t)})$ must be also in $V$ which implies that this
vector cannot be general. A contradiction.

To assure  ($\#'$)  one can choose $x_1, \ldots , x_{n-1}$ so that $x_i=x(it_0)$ for some nonzero $t_0$ and choose an automorphism $\alpha$ so that
after its application the vectors ${\rm d}\rho_0 (\delta_{1,x_i})-{\rm d}\rho_0 (\delta_{1,x_{i-1}}), \,i=1, \ldots , n-1$ form a basis of $T_{q_0}Q_0$
(where $x_0=x$).

(b) Moreover, the proof of Proposition \ref{fm.p1} implies that one can choose $\delta_1$ so that condition ($\#'$) holds not only for the orbit $O_1$ of one general point $x \in X$
but simultaneously for the orbits of any finite set of general points in $X$.

\erem

We need some further facts from \cite{FKZ-GW}.

\bnota\label{fm.n1} (a) Denote by $U^i$ the unipotent one-parameter subgroup associated with $\delta_i$ from Proposition \ref{fm.p1} and
for every $f \in \Ker \delta_0 \setminus \Ker \delta_1$ (resp.  $g \in \Ker \delta_1 \setminus \Ker \delta_0$) denote by $U_f^0$ (resp. $U_g^1$)
the one-parameter group associated with the replica $f\delta_0$ (resp. $g\delta_1$).

(b) To any sequence of invariant functions
\be\label{seq}
\cF=\{f_1,\ldots,f_s, g_1, \ldots,g_s\},\,\,\,\,\mbox{where}
\,\,\,\, f_i\in\Ker\de_1\backslash\Ker
\de_0\,\,\,\,\mbox{and}\,\,\,\, g_i\in\Ker\de_0\backslash\Ker
\de_1\,,
\ee
we associate an algebraic family
of automorphisms $\C^{2s}\to\Aut(X)$ defined by the product
\be\label{00121}
U^\cF=U^1_{f_s}\cdot
U^0_{g_s}\cdot\ldots\cdot U^1_{f_1}\cdot
U^0_{g_1}\subseteq H\,.
\ee
More generally, given a tuple
$\kappa=(k_i,l_i)_{i=1,\ldots,s}\in\N^{2s}$ the product
\be\label{001210}
U_\kappa=U_\kappa^\cF=
U^1_{f_s^{k_s}}\cdot
U^0_{g_s^{l_s}}\cdot\ldots\cdot U^1_{f_1^{k_1}}\cdot
U^0_{g_1^{l_1}}
\subseteq H\,
\ee
yields as well an algebraic family of automorphisms.
\enota

\bprop\label{fm.p2} {\rm (\cite[Corollary 4.4]{FKZ-GW})} There is a finite collection of invariant functions $\cF$ as in (\ref{seq})
such that for any sequence
$\kappa=(k_i,l_i)_{i=1,\ldots,s}\in\N^{2s}$ the algebraic family of automorphisms
$U_\kappa$ as in (\ref{001210})
has a dense open orbit in $X$.
This orbit $O(U_\kappa)$ coincides with $O(H)$ and so
does not depend on the choice of $\kappa\in\N^{2s}$.\eprop

\brem\label{fm.r2} Let 
$\tilde \cF=\{\tilde f_1, \ldots, \tilde f_r, f_1,\ldots,f_s, \tilde g_1, \ldots, \tilde g_r, g_1, \ldots,g_s\}$ where $f_i$ and $g_i$ are as in Notation \ref{fm.n1}
while $\tilde f_i\in\Ker\de_1\backslash\Ker
\de_0$ and $\tilde g_i\in\Ker\de_0\backslash\Ker
\de_1$. Note that if $O(U_\kappa^\cF)=O(H)$ as in Proposition \ref{fm.p2} then one has also  $O(U_\kappa^{\tilde \cF})=O(H)$.

\erem

\bnota\label{fm.n2}
(a) Given a one-parameter group $U\in \cU(X)$ we let $U^*=U\backslash\{\id\}$.
Given a collection $\cF=\{f_1,\ldots,f_s, g_1, \ldots,g_s\}$ of invariant functions as in Notation \ref{fm.n1}
and $U_\kappa=
U^1_{f_s^{k_s}}\cdot
U^0_{g_s^{l_s}}\cdot\ldots\cdot U^1_{f_1^{k_1}}\cdot
U^0_{g_1^{l_1}}$ as in \eqref{001210},  we let
$$
U_\kappa^*=
U^{1*}_{f_s^{k_s}}\cdot
U^{0*}_{g_s^{l_s}}\cdot\ldots\cdot U^{1*}_{f_1^{k_1}}\cdot
U^{0*}_{g_1^{l_1}}\,.
$$
(b) Consider $\delta_0$ as in Notation \ref{fm.n1} and its partial quotient morphism $\rho_0 : X \to Q_0$. Then it can be extended to a proper morphism $\bar \rho_0: \bX \to Q_0$. There is exactly one 
(so-called horizontal) irreducible component $D_0$ of  the variety  $\bX \setminus X$ for which the restriction of the morphism $\bar \rho|_{D_0} : D_0 \to Q$ is birational  \cite{FKZ-GW}.
\enota

\bprop\label{fm.p3} Let the assumption of Proposition \ref{fm.p2} hold and $Z$ be a closed subvariety of $X$ of codimension at least 2. 
Then the integers $k_1, l_1, k_2, l_2, \ldots , k_s, l_s$ in Notation \ref{fm.n2}
can be chosen so that there exists a proper subvariety $R\subset D_0$ such that for every element $\alpha \in U_\kappa^*$ 
the closure $\bZ_\alpha$ of $\alpha (Z)$ in $\bX$
meets $D_0$ from Notation \ref{fm.n2} along $R$ only. 
\eprop

\bproof The statement is a special case of \cite[Proposition 4.11]{FKZ-GW}). 
\eproof

\blem\label{fm.l2} Let the assumption of Theorem \ref{fm.t1} hold  and let $U_\kappa, \rho_0, \bar \rho_0, D_0$ be as Proposition \ref{fm.p3} with $\delta_0=\sigma$.
Then there exists an automorphism $\alpha \in \Aut (X)$ such that

{\rm (a)}  the conclusions of Lemma \ref{fm.l1} are satisfied for every point 
$z_i'=\alpha(z_i)$ in $Z'=\alpha (Z)$ and $\rho_0 (z_i') \notin  \bar \rho_0 (\overline {Z'} \cap D_0) \subset  \bar \rho_0 (R')$ 
where $R'$ plays the same role for $Z'$ as $R$ for $Z$ in Proposition \ref{fm.p3}.

Furthermore, let $O_i=\rho_0^{-1}(\rho_0 (z_i'))\, \, i=1, \ldots , m$, $M_i \subset O_i$ be a finite subset, and $B_y$ be the analytic branch of $Z'$ at $y \in M_i$.
Then for every $i$ and every such a point $y$

(b)  the vector  ${\rm d}\rho_0 (\delta_{1,z_i'})-{\rm d}\rho_0 (\delta_{1,y})$ is not contained in the tangent cone of $\rho_0(B_y)$ at $\rho_0 (z_i')$.
\elem

\bproof Choose an automorphism $\beta \in \Aut (X)$ such that the conclusions of Lemma \ref{fm.l1} hold for
every point $z_i''=\beta (z_i)$ of the variety $Z''=\beta (Z)$. That is,  each $z_i''$ is a general point of $X$, 
$\rho|_{Z''} : Z'' \to \rho (Z'')$ is a local embedding at $z_i''$,
and $\rho (z_i'')$ is contained in a smooth analytic branch of $\rho (Z'')$
which is in turn contained in the smooth part of $Q$. Furthermore, perturbing $\alpha$ by the virtue of \cite[Proposition 4.14 and Remark 4.16]{AFKKZ}
we can suppose that condition (b) holds with $Z'$ and $z_i'$ replaced by $Z''$ and $z_i''$.
Note that this conditions are preserved by perturbation of $Z''$ by any
automorphism $\gamma \in U_\kappa^*$ sufficiently close to the identical map (in the compact-open topology).
That is the conclusions of Lemma \ref{fm.l1} are still valid for $Z'=\gamma (Z'')$ and $z_i'=\gamma (z_i'')$.
Note also that by Proposition \ref{fm.p3} $\rho_0 (z_i') \notin  \bar \rho_0 (R')$ for any general $\gamma$.   Thus one has now (a) since the orbit of any general point
$z_i''$ under the action of  $U_\kappa^*$ is open by Proposition \ref{fm.p2}. Since condition (b) is preserved for general $\gamma$ we are done.

\eproof


\subsection{Proof of Theorem \ref{fm.t1}.} 
Treat $\sigma, z_1\in Z$, and $ \rho^{-1}(\rho(z_1))$ as $\delta_0, z_1'\in Z'$ and  $O_1$ in Lemma \ref{fm.l2}. 
Choose $\delta_1$ so that the assumptions of Proposition \ref{fm.p1} and Remark \ref{fm.r1} are satisfied for a general point $y \in O_1\simeq \C$, and, 
in particular,  the vector ${\rm d}\rho_0 (\delta_{1,x(t)})-{\rm d}\rho_0 (\delta_{1,y(t)})$ is general 
in $T_{q_0}Q_0$ for general moment of time $t$. Suppose that a finite set $M\subset O_1$ consists of all non-general points. 

Choose $U_\kappa, R'$, and $Z'$ as in the proof of Lemma \ref{fm.l2} (in particular condition (b) is valid for points in $M$). Keep in mind that for general $\alpha \in U_\kappa$ the conclusions of Lemma \ref{fm.l1}
are satisfied. That is, $z_1'$ is a general point of $X$, 
$\rho_0|_{Z'} : Z' \to \rho_0 (Z')$ is a local embedding at $z_1'$,
and $\rho_0 (z_1')$ is contained in a smooth analytic branch of $\rho_0 (Z')$
which is in turn contained in the smooth part of $Q$. For Theorem \ref{fm.t1} we only need to establish that 
$\rho_0|_{Z'}$ is a birational morphism on its image and $\rho_0 (Z')$ is unibranch at $q_1=\rho_0(z_1')$ 
which is now equivalent to the fact that for the set  $\rho_0^{-1}(q_1) \cap Z'$ contains no other points but $z_1'$. 

We are going to show that this equivalent fact is true under replacement of $Z'$ by $\alpha (Z')$ for a general  $\alpha \in U_\kappa$.
Suppose that $f_1=g_1=1$ for $f_1,g_1 \in \cF$ from Notation \ref{fm.n1} which is
allowed by Remark \ref{fm.r2}. The let us  look for this $\alpha$ in the form $\alpha =\beta \alpha_1\alpha_0$ where $\alpha_0 \in U^{0*}_{g_1^{l_1}}$
and  $\alpha_1 \in U^{1*}_{f_1^{k_1}}$. Furthermore, we suppose that
$\beta$ is as close to the identical automorphism as we wish and therefore $\alpha (Z')$ is as close to $\alpha_1 \alpha_0 (Z')$ near point $z_1'$ 
as we wish. 

Assume that exists  a point $y_1\in \rho_0^{-1}(q_1) \cap (Z'\setminus z_1')$ (since otherwise we are done). 
Let $B$ be the analytic branch of $Z'$ at $y_1$ and $B_0=\rho_0 (B)$ be the analytic branch of $\rho_0 (Z')$ at $q_1$. Note that 
for any element $\alpha_0$ we have equality $B_0=\rho_0 (\alpha_0 (B))$. 

{\em Case 1}. Suppose that $y_1$ does not belong to $M$. By Remark \ref{fm.r1} after application of a general $\alpha_0$
the vector $$\nu ={\rm d} \rho_0 ({\rm d} \alpha_1 (\delta_{1,y_1}))- {\rm d} \rho_0( {\rm d} \alpha_1 (\delta_{1,z_1'}))$$ is general. 
That is, $\nu$  is not contained in the tangent cone $C_{q_1}B_0$ to $B_0$ at $q_1$. 

Denote by $B_0'$ the image of $\alpha (B)$ under $\rho_0$. Observe that if $\alpha_1$ is close to the identity isomorphism then up to infinitesimals  $\rho_0(\alpha_1\alpha_0 (z_1'))$ 
 changes from $q_1$ in direction vector $c{\rm d} \rho_0( {\rm d} \alpha_1 (\delta_{1,z_1'}))$ (where $c$ is a nonzero coefficient) 
 while the cone $C_{\rho_0(\alpha_1(y_1))}B_0'$
 is obtained from the cone $C_{q_1}B_0$ via the shift by vector $c{\rm d} \rho_0( {\rm d} \alpha_0 (\delta_{1,y_1}))$. 
 
This implies that for a general $\alpha_1$
the point $\rho_0 ( \alpha_1 \circ \alpha_0 (z_1'))$ is not contained in $\rho_0 ( \alpha_1 \circ \alpha_0 (B))$, i.e. $ \alpha_1 \circ \alpha_0 (B)$ 
does not meet $\rho_0^{-1}( \rho_0 (\alpha_1 \circ \alpha_0(z_1'))$. 
Since $\beta$ is close to the identical map as we wish the same is true for general $\alpha =\beta \circ \alpha_1 \circ \alpha_0$,
i.e.  $ \alpha (B)$ 
does not meet $\rho_0^{-1}( \rho_0 (\alpha (z_1')))$. Therefore, $\alpha (Z')$ does not meet $\rho_0^{-1}( \rho_0 (\alpha (z_1')))$ near any
point $y_1\in \rho_0^{-1}(q_1) \cap (Z\setminus z_1)$ for general $\alpha$. 

{\em Case 2}. Suppose that $y_1\in M$ is not general in $\rho_0^{-1}(q_1)$. Then we can suppose that $\nu$  is not contained in the tangent cone $C_{q_1}B_0$ to $B_0$ at $q_1$
by Lemma \ref{fm.l2} (b) which yields the same conclusion that $\alpha (Z')$ does not meet $\rho_0^{-1}( \rho_0 (\alpha (z_1')))$ near any
point $y_1\in \rho_0^{-1}(q_1) \cap (Z\setminus z_1)$ for general $\alpha$. 

Note also that since $\rho_0 (z_1') \notin  \bar \rho_0 (\overline {Z'} \cap D_0) \subset  \bar \rho_0 (R')$ the points in
$\rho_0^{-1} (\rho_0 (\alpha (z_1')) \cap \alpha (Z')$ depend upper-semicontinuously on $\alpha$.
Hence the previous argument shows that $\rho_0^{-1} (\rho_0 (\alpha (z_1')) \cap \alpha (Z' \setminus z_1')=\emptyset$, i.e. we have 
the statements (i)-(iv) of Theorem \ref{fm.t1} for $z_1$ and general $\alpha$.

Note that in this argument we used only the fact that the restriction of $\delta_1$ to the orbit $O_1$ of $z_1$ under $U^0$ satisfies condition ($\#'$) from Remark \ref{fm.r1}.
However, the same Remark \ref{fm.r1} shows that $\delta_1$ can be chosen so that it satisfies ($\#'$) for the orbit of every $z_i, \, i=1, \ldots , m$. Thus, for a general $\alpha$
we have the statements (i)-(iv) of Theorem \ref{fm.t1}  for each $z_i$ which concludes the proof.
$\square$

\section{First facts about sprays on flexible manifolds.}

\bthm\label{ff.t1} Every flexible manifold $X$ is elliptic.
\ethm

\bproof  Recall that for every $x \in X$ and each nonzero vector $v \in T_xX$ there is a locally nilpotent vector field $\sigma$ on $X$
for which the value $\sigma_x$ of $\sigma$ at $x$ coincides with $v$ \cite[Corollary 4.3]{AFKKZ}. In particular, we find locally nilpotent vector fields $\sigma_1, \ldots , \sigma_n$ for
which $\sigma_{1,x}, \ldots , \sigma_{n,x}$ is a basis in $T_xX$. This implies that there is an open Zariski dense subset $U$ of $X$ so that
for every point $y \in U$ the vectors $\sigma_{1,y}, \ldots , \sigma_{n,y}$ is a basis in $T_yX$.  Note that $\dim (X \setminus U)=m\leq n-1$.
Finding locally nilpotent vector fields that form a basis at general points of each of the component of $X \setminus U$ we can extend
our sequence of vector fields to $\sigma_1, \ldots , \sigma_n, \sigma_{n+1}, \ldots , \sigma_l$ such that there is a Zariski open set $V\supset U$
for which $\dim (X \setminus V)<m$ and such that for every $y \in V$ these fields generate $T_yX$. Thus using induction by dimension
we can suppose that $\sigma_1, \ldots , \sigma_n, \sigma_{n+1}, \ldots , \sigma_l$ generate $T_yX$ at every $y \in X$.

Let $U^i$ be the one-parameter group of algebraic automorphism associated with $\sigma_i$ and let $U^i_t$ be the element of this group
for the value of the time parameter $t$. Consider the trivial vector bundle $\pi : E \to X$ of rank $l$, i.e. for every $x \in X$ the fiber $E_x=\pi^{-1}(x)$ is isomorphic to $\C^l$ with
coordinates $(t_1, \ldots , t_l)$.

Define the morphism $s: E \to X$ by the formula
$$(x,t_1, \ldots , t_n)\to U_{t_1}^1\circ \cdots \circ U_{t_l}^l (x).$$
By construction, $s$ is a dominating spray and we are done.

\eproof 

By Corollary \ref{sas.c1} we have the following.

\bcor\label{ff.c1} Every locally flexible algebraic manifold is subelliptic.

\ecor

\bnota\label{ff.n1} Let $\pi : \tX \to X$ be the blowing up of an affine manifold $X$ along a closed smooth algebraic submanifold $Z$ of codimension $k+1 \geq 2$.
Suppose that $E$ is the exceptional divisor of $\pi$, i.e. $\pi|_E : E \to Z$ is a locally trivial fibration with fiber $\PP^k$.
Note that $\C [\tX ]=\C [X]$ and $\tX$ is a semi-affine manifold. In particular, the notion of locally nilpotent derivation (= vector field) on $\tX$ is well-defined.
\enota

\bprop\label{ff.p1} For every $z \in Z$, $x \in \pi^{-1}(z)$, and a nonzero vector $w \in T_x\tX$ tangent to $\pi^{-1}(z)$ there exists a locally nilpotent
vector field $\delta$ on $\tX$ for which $\delta_x =w$. Furthermore, the subgroup of automorphisms of $\tX$ preserving $\pi^{-1}(z)$ acts transitively
on  $\pi^{-1}(z)$.

\eprop

\bproof Let $\rho : X \to Q$ be a partial quotient associated with a nonzero locally nilpotent $\sigma$, $x$ be general point in $X$ and $q=\rho (x)$ a general point in $Q$.
In particular these points are smooth and one can choose local analytic  coordinate systems at them. 
By \cite[Theorem 4.14 and Remark 4.16]{AFKKZ} one can also choose an automorphism $\alpha$ which sends $z $ to $x$
such that in  local analytic coordinate systems $(v_{0}, \ldots , v_{n-1})$ (resp. $(u_{0}, \ldots , u_{n})$)  at $q\in Q$ (resp. $x \in X$) one has $\alpha (Z)$ given by $u_{n-k}=\ldots = u_{n}=0$
and $\rho^* (v_{j})=u_{j}$ for $j\leq n-1$.  We replace $Z$ by $\alpha (Z)$ and $z$ by $\alpha (z)$ to make the argument local.
In particular, $\pi^{-1}(z)\simeq \PP^k$ has homogeneous coordinates $U_{n-k}:U_{n-k+1}: \ldots : U_n$ such that $u_iU_j=u_jU_i$ for $n-k \leq i,j \leq n$. 
Without loss of generality consider the case when 
the vector $w$ in $T_{\pi^{-1}(z)}$ is tangent to the line $L$ with fixed relation
$U_{n-k+1}: \ldots : U_{n-1}:U_n$ where $U_n \ne 0$ and arbitrary $U_{n-k}$. 
Note that at the origin $x_0$ of the local coordinate system $\sigma_{x_0}$ is proportional
to the vector $\p / \p u_{n}$, i.e. we can suppose that  $\sigma_{x_0}=\p / \p u_{n}$. Since $u_{n-k}=\rho^*(v_{n-k})\in \Ker \sigma$ we see that $ u_{n-k}\sigma$ is also locally nilpotent.
Denote by $\Phi$ the automorphism $\Phi=\exp (tu_{n-k}\sigma)$ of $X$ for some value of parameter $t \in \C$. By  \cite[Lemma 4.1]{AFKKZ} we have
\be\label{ff.eq1} {\rm d}_{x_0} \Phi (\nu) =\nu +t {\rm d}u_{n-k}(\nu) \p / \p u_{n}\ee for every $\nu \in T_{x_0}X$. Since $u_{n-k}\sigma$ vanishes on $Z$ it can be lifted as
a locally nilpotent derivation $\delta$ on $\tX$. Furthermore, Formula (\ref{ff.eq1}) shows that the elements of the flow of $\delta$ preserve $\pi^{-1}(z)\simeq \PP^k$
and act on as elementary transformations of form   $(U_{n-k}:U_{n-k+1}: \ldots : U_n ) \to ((U_{n-k}+tU_n):U_{n-k+1}: \ldots : U_n)$.  
That is,  the action induced by $\delta$ is a translation along the affine line $\C \simeq L\setminus \{ U_{n-k}=\infty \}$
which yields the first statement. The fact that elementary transformations generate a special linear group implies the second statement and we are done.

\eproof

\section{Affine modifications}

The next definition of affine modifications and their properties can be found in \cite{KZ}.

\bdefi\label{am.d1}
A birational morphism $\varphi : X' \to X$ of affine algebraic varieties is called  an affine modification. In particular, one can find effective reduced
divisors $D\subset X$ and $E \subset X'$ such that the restriction $\varphi|_{X' \setminus E} : X' \setminus E \to X \setminus D$ is an
isomorphism. Though these divisors are not determined uniquely it will be clear from the context below what they are.
We call $E$ the exceptional divisor of the modification and $D$ the divisor of the modification. Furthermore, we consider only
the cases when $D$ is principal, i.e. $D=f^*(0)$ for some regular function $f$ from the ring $A=\C [X]$ of regular functions on $X$.
In this case the ring $A'=\C [X']$ of regular functions on $X'$ is generated over $A$ by functions of form $g/f$ where $g$ runs over
an ideal $I$ of $A$, i.e. $A'=A[I/f]$. Even for a fixed $f$ this ideal may not be unique and to remove this ambiguity we suppose that
it is the largest ideal for which  $A'=A[I/f]$. This largest ideal will be called the ideal of the modification. The center of modification
(for fixed $D$ and $E$) is the closure $Z$ of $\varphi (E)$ in $X$.
\edefi

\brem\label{am.r1}

 The geometrical meaning of the modification is the following.
One consider the blowing up $\pi : \tX \to X$ of $X$ along the ideal $I$ and obtain $X'$ by removing from $\tX$ those divisors
on which the zero multiplicity of $f$ is greater than the zero multiplicity of at least one function $g$ from $I$ (and letting $\varphi = \pi|_{X'}$). In particular,
if $D$ and $Z$ are smooth and the ideal of the modification coincides with the defining ideal $I(Z)$ of the center then
$X'$ is the complement in $\tX$ to the proper transform of $D$.

\erem

\blem\label{am.l1} Let $\rho : X \to Q$ be a dominant morphism of normal irreducible affine algebraic varieties, $D$ be an effective principal divisor 
in $Q$, $Q_0=Q\setminus D$, and $X_0=\rho^{-1}(Q_0)$. Suppose that $X_0$ is isomorphic to $Q_0 \times \C^m$ over $Q_0$.
Then 

{\rm (a)} there exists an affine modification $\varphi : X \to Q \times \C^m$ whose restriction over $Q_0$ is an isomorphism;

{\rm (b)} for every locally nilpotent vector field $\sigma$ on $X_0$ tangent to the fibers of $\rho|_{X_0}$ there exists an equivalent\footnote{We call two locally
nilpotent vector fields equivalent if as derivations they have the same kernel.}
locally nilpotent vector field $\delta$ that extends regularly to $X$;

{\rm (c)} furthermore, for every $q \in Q_0$ the restrictions of the fields $\sigma$ and $\delta$ to the fiber $\rho^{-1}(q)$ differ by 
a nonzero constant factor.  

\elem

\bproof Note that the isomorphism $X_0 \simeq Q_0 \times \C ^m$ has a coordinate form $(\rho , h_1, \ldots , h_m)$ where each $h_i$ is a regular function
on $X_0$. If these functions extend regularly to $X$ then it suffices to put $\varphi = (\rho , h_1, \ldots , h_m)$.
Otherwise, consider a regular function $g \in \C [Q]$
for which $D =g^*(0)$. Note that the extension of $h_i$ to $X$ can have poles only on the divisor $\rho^{-1}(D)$. Thus for sufficiently
large $k_i$ the function $g^{k_i}h_i$ is regular on $X$ (because of normality). Replacing every $h_i$ by $g^{k_i}h_i$ we see that the same $\varphi$ yields 
the desired affine modification in (a).

Similarly, since by the assumption $g\in \Ker \sigma$ the field $\delta =g^k \sigma$ is also locally nilpotent and for $k$ large enough
it extends regularly to $X$. Thus we have (b). For (c) it suffices to observe that $g$ does not vanish on $Q_0$. 

\eproof

\bprop\label{am.p1}
Let the assumption of Lemma \ref{am.l1} hold, $D=g^{-1}(0)$, and $\tau : Q\times \C_{\bar u}^m \to Q \times \C_{\bar u}^m$ 
(where $\bar u =(u_1, \ldots , u_m)$ is a coordinate system on $\C^m$) be a birational morphism over $Q$ such that
$\tau^*(\bar u)=\bar u +g^k\bar e$ where $\bar e$ is a regular function on $Q$ with values in $\C^m$. 
Then for $k$ large enough this endomorphism can be lifted to
a birational morphism $\theta : X \to X$, i.e. the following commutative diagram holds

\[  \begin{array}{ccc} X&  \stackrel{{\theta}}{\rightarrow} & X \\
\, \, \, \, \downarrow {\varphi}  & & \, \, \, \, \downarrow { \varphi} \\
Q\times \C^m&  \stackrel{{\tau}}{\rightarrow} & Q \times \C^m.

\end{array} \]
Furthermore, $\theta$ maps $\varphi^{-1}(D)$ isomorphically on itself.

\eprop

\bproof By Lemma \ref{am.l1} one can consider an affine modification $\varphi : X \to Q \times \C^m$, i.e. $\C [X]$ is generated over $\C [Q \times \C^m]$ by
elements of form $I/g^k$ where $k$ is some natural number and $I$ is an ideal in $\C [Q \times \C^m]$ generated by $g^{k_0}$ and 
elements $g_1, \ldots, g_k\in\C[Q\times \C^m]$.
By the assumption $\tau^*$ transfer $g$ to $g$ and $I$ into the ideal $J$ generated by $g^{k_0}$ and $\tau^*(g_1), \ldots , \tau (g_k)$.
Choose $k \geq k_0$. By Taylor expansion one has $\tau (g_i)= g_i +g^kh_i$ where $h_i$ is regular on $Q\times \C^m$.
Hence, $J$ is contained in $I$. Now the first statement follows from \cite[Proposition 2.1]{KZ}.

Note also that $\tau$ is invertible in an \'etale neighborhood of $D$ in $Q \times \C^m$. Hence $\theta$ is invertible
in an \'etale neighborhood of $\varphi^{-1}(D)$ in $X$ which yields the second statement.
\eproof

\bcor\label{am.c1} Let $Q$ and $X$ be algebraic varieties, $X_0=Q\times \C_{\bar u}^m$, $\varphi : X\to X_0$ be an affine modification
over $Q$ with divisor $D \subset X_0$ given by the zeros of a regular function $g \in \C[Q] \subset \C [X_0]$, and $f \in \C [Q]\subset \C [X_0]$ 
be another regular function that has disjoint zeros with those of $g$. Furthermore, suppose that $f|_D\equiv 1$ with multiplicity $k$. 
Let $\tau : X_0 \to X_0$ be the birational morphism over $Q$ for which $\tau^* (\bar u)=f(q) \bar u $. Then for $k$ large enough there exists a modification $\theta : X \to X$
over $Q$ such that the commutative diagram and the second statement from Proposition \ref{am.p1} hold.

\ecor

Now we have the following observation. 

\bthm\label{am.t1} Let the assumptions of Corollary \ref{am.c1} hold, $Z=\varphi^{-1}(Z_0)$ where $Z_0$ is given
in $X_0$ by the equations $\bar u=\bar 0$ and $f=0$. Suppose that $\pi : \tX \to X$ is the blowing of $X$ up along the center $Z$.
Then $\tX$ contains a Zariski open set $X'$ isomorphic to $X$ such that $X'\cap E$ is dense in the exceptional
divisor $E=\pi^{-1}(Z)$ of $\pi$.

\ethm

\bproof Consider the affine modification $\psi : X' \to X$ along the divisor $f^{-1}(0)$ with center $Z$. That is, $X'$ is
a Zariski open subset of $\tX$ with $X'\cap E$ dense in $E$. On the other hand $\theta : X \to X$ is also an affine modification
with center $Z$. Its exceptional divisor contains $f^{-1}(0)$ and by Corollary \ref{am.c1} it contains nothing else. This yields the
desired isomorphism $X\simeq X'$.

\eproof

\bnota\label{am.n2} Let $X$ be a quasi-affine algebraic manifold, $Z$ be a submanifold of $X$ which is a strict complete
intersection. That is, the defining ideal $I$ of $Z$ is generated by regular functions $f:=g_0,g_1, \ldots , g_k$ where $k+1$ is
the codimension of $Z$ in $X$. For $l<k$ consider the strict complete intersection $Z_1 \subset X$ given by $f=g_1=\ldots = g_l=0$.
Suppose that $\varphi : X_1' \to X$ is the modification with center $Z_1$ and divisor $D=f^*(0)$
while $\pi_1: \tX_1 \to X$ (resp. $\pi : \tX \to X$) is the blowing up of $X$ with center at $Z_1$ (resp. $Z$), i.e. $X_1'$ can be viewed as a Zariski open subset of $\tX_1$.
Suppose that $\tZ_1\subset X_1'$ is the complete strict intersection given by $f\circ \varphi =g_{l+1} \circ \varphi = \ldots = g_k \circ \varphi =0$
and $\pi' : \tX_1' \to X'$ is the blowing up of $X_1'$ with center at $\tZ$.

\enota

\bprop\label{am.p2} There is a natural birational morphism $\theta : \tX_1' \to \tX$ such that $\theta (\tX_1')$ meets the exceptional divisor $E$ of $\pi$
along a Zariski  dense open subset of $E$.

\eprop

\bproof Recall that $\tX$ can be viewed as the submanifold of $X \times \PP^k$ given by equations $U_ig_j=U_jg_i$ for $i,j=0, \ldots , k$
where $(U_0: U_1: \cdots : U_k)$ is a homogeneous coordinate system on $\PP^k$.
Similarly $\tX_1$ is the submanifold of $X \times \PP^l$ given by equations $V_ig_j=V_jg_i$ for $i,j=0, \ldots , l$
where $(V_0: V_1: \cdots : V_l)$ is a homogeneous coordinate system on $\PP^l$. Note that $X_1'$ is given in $\tX_1$
by the equation $V_0=1$ and  $\tX_1'$ is the submanifold of $X_1' \times \PP^{k-1}$ given by equations $W_ig_j=W_jg_i$ for $i,j=0, l+1, l+2\ldots , k$
where $(W_0: W_{l+1}: \cdots : W_k)$ is a homogeneous coordinate system on $\PP^{k-l}$.

Note that one has the natural birational map $\theta : \tX_1'  \rDashto \tX$ over $X$ which is an isomorphism over $X\setminus Z_1$ and regular over $X \setminus Z$.
Thus one needs to check only the regularity over $Z$. However, one can see that over $Z$ this map is automatically given by the following  
\[ \begin{array}{c}  [(V_0: V_1: \cdots : V_l), (W_0: W_{l+1}: \cdots : W_k)] \to (U_0: U_1: \cdots : U_k)=\\
(V_0W_0 : V_1W_0: \cdots : V_lW_0 : W_{l+1}V_0 : \cdots : W_kV_0)\end{array}\]
which is regular since $V_0=1$. Note also that when $W_0=1$ then this morphism is an embedding. This yields the desired conclusion
about the density of the intersection of $E$ and $\theta (\tX_1')$.

\eproof

\section{$k$-Flexibility}

\bdefi\label{kf.d1} A flexible quasi-affine manifold $X$ will be called $k$-flexible for $k>0$ if there exists a 
morphism $\rho : X \to Q$ into a normal affine algebraic variety $Q$ such that over a Zariski open dense subset $Q_0$ of $Q$ the variety
$\rho^{-1}(Q_0)$ is 
 isomorphic to $Q_0 \times \C^k$ over $Q_0$.  

\edefi

\brem\label{kf.r1} (i) Note that for $l>k$ each $l$-flexible variety is automatically $k$-flexible.

(ii) Every flexible manifold $X$ is, of course, 1-flexible, since one can consider any partial quotient morphism
of $\rho : X \to Q$. Then for some $Q_0$ as above $\rho^{-1}(Q_0)\to Q_0$ is a locally trivial $\C$-fibration. Requiring that
$Q_0$ is affine one can guarantee that it is in fact an line bundle. Hence  removing from $Q_0$ a divisor we make 
$\rho^{-1}(Q_0)$ the desired direct product.

(iii) For $k=2$ it is also enough to require that $\rho : X \to Q$ has general fibers isomorphic to $\C^2$. The existence
of a desired $Q_0$ follows from \cite{KZ01}.

\erem

\bexa\label{kf.e1}  Consider a hypersurface $H$ given by $uv =p(\bar x)$ in $\C^{n+2}$ where $u, v$, and $\bar x=(x_1, \ldots , x_n)$ are coordinates
on $\C^{n+2}$. If the zero locus of $p$ in $\C^n_{\bar x}$ is smooth then $H$ is a flexible manifold \cite{KZ}. 
Note that $\dim H=n+1$ and $H$ is $n$-flexible.
\eexa

\bprop\label{kf.p1} Let $X=SL (n)$, i.e. $\dim X=n^2-1$. 
Then $X$ is $k$-flexible where $k=n^2-n$. 

\eprop

\bproof Let $A=[a_{ij}]_{i,j=1}^n$ be a matrix from $SL(n)$
and let $\{ A_{ij} \}$ be cofactors of this matrix.  Consider the natural morphism of $\pi :X\to H$
into the hypersurface $H$ given by the equation 
\be\label{kf.eq1} a_{11}A_{11}+a_{12}A_{12} + \ldots + a_{1n}A_{1n}=1\ee
in $\C^{2n}$ with coordinates $(a_{11}, a_{12}, \ldots , a_{1n}, A_{11}, A_{12}, \ldots , A_{1n})$. Note that $\dim H=2n-1$.
Let $A'$ (resp. $A''$) be the matrix obtained from $A$ by removing the first row (resp. the first row and the first column). 
Consider the action of $SL(n-1)$ on $X$ such that
for $B \in SL(n-1)$ the matrix $B.A$ is obtained by replacing $A'$ by $BA'$
while keeping the first row in $A$ intact. This action is free and it preserves the fibers of $\pi$ which are therefore of dimension at least $(n-1)^2-1$.
Observing  the equality $n^2-1=(n-1)^2-1 +(2n-1)$ one can see now that the fibers of $\pi$ are nothing but the orbits of this action.

Furthermore, let $H_0$ be the complement to the zero locus of $A_{11}$ in $H$. Since for every point $A$ in  $X':= \pi^{-1}(H_0)$ the determinant of $A''$ is nonzero, dividing the first row of $A''$
by this determinant we see that $X'$ is 
 isomorphic to $H_0 \times SL(n-1)$. Note also that by Formula (\ref{kf.eq1}) $H_0\simeq \C^* \times \C^{2n-2}$ where 
 $\C^* \subset \C$ is equipped with the coordinate $A_{11}$. Hence $X' \simeq  \C^* \times SL(n-1) \times \C^{2n-2}$.

One can check that $SL(2)$ is 2-flexible.
Thus by the induction assumption we can suppose that there are a morphism $\rho_{n-1}: SL(n-1)\to Q_{n-1}$ into a normal affine algebraic variety
$Q_{n-1}$ and a Zariski dense open subset $Q_{n-1}^0$ in it for which the Zariski open subset $\rho^{-1}_{n-1} (Q_{n-1}^0)$ of  
$SL(n-1)$ is 
 isomorphic to $Q_{n-1}^0\times \C^l$ with $l=(n-1)^2-n+1$. This yields a Zariski open subset $X_0$ of $X'$ 
  isomorphic to
 $Q_n^0\times \C^k$ where $k=l+2n-2=n^2-n$ and   $Q_n^0:=\C^*\times Q_{n-1}^0$. 
 Hence it suffices to modify this isomorphism  $\psi : X_0\simeq Q_n^0\times \C^k$ so that it becomes a restriction of an affine modification $X\to Q_n\times \C^k$ over some variety $Q_n$
 containing $Q_n^0$ as a Zariski open subset. 
  
 Consider morphism $\tau =(A_{11}, \rho_{n-1}) : X' \to P:=\C\times Q_{n-1}$.  Without loss of generality we can suppose that  $P\setminus Q_n^0$ is the principal divisor in $P$. 
 Then modifying $\psi$ we can extend it to a morphism  $\varphi : X' \to P \times \C^k$  over $P$ by Lemma \ref{am.l1}(a).
 
 Treat now $Q_{n-1}$ as a closed subvariety of $\C^m$, i.e. $P$ is a closed subvariety of $\C^{m+1}$ with coordinates $(u_0,u_1, \ldots , u_{k})$
 where the lift of $u_0$ to $X'$ coincides  with $A_{11}$.  
 Consider composition of $\tau$ with endomorphism $\theta: \C^{m+1} \to \C^{m+1}$ given by $$(u_0,u_1, \ldots , u_{k})\to (u_0,u_0^Nu_1, \ldots , u_0^Nu_{k})$$
 where $N$ is natural. Since the restriction of $\theta$ to the complement to $u_0=0$ is an isomorphism we see that $\tau (Q_n^0)$ is isomorphic
 to its image $\theta\circ \tau (Q_n^0)$ which is, therefore, isomorphic to $Q_n^0$. Furthermore, since  $X\setminus X'$ is the zero locus of function $A_{11}$ on $X$,
 for sufficiently large $N$ this morphism $\theta \circ \tau$ extends to a morphism $\rho_n: X \to Q_n$ where $Q_n$ is the closure of $\theta\circ \tau (Q_n^0)$.
 
 Applying Lemma \ref{am.l1}(a) again we obtain the desired change of isomorphism $\psi : X_0\simeq Q_n^0\times \C^k$ so that it because an affine
 modification $X \to Q_n\times \C^k$ over $Q_n$ which concludes induction and the proof.

\eproof

\bthm\label{kf.t1} Let $X$ be a $k$-flexible quasi-affine manifold and let $\rho : X \to Q$ and $Q_0$ be as in Definition \ref{kf.d1}.
Suppose that $D_0=Q \setminus Q_0$ is a principal divisor,
$Z$ is a closed submanifold of $X$ of codimension $k+1$
such that $\rho|_Z : Z \setminus \rho^{-1}(D_0)\to Z_0$ is an isomorphism for $Z_0 =\overline{\rho (Z)}\setminus D_0$ which is a principal divisor in $Q_0$.
Let $\pi : \tX \to X$ be the blowing of $X$ up along the center $Z$.

Then 

{\rm (i)} $\tX$ contains a Zariski open set $X_0'$ 
 isomorphic to $X_0=Q_0 \times \C^k$ 
such that $X_0'\cap E$ is dense in the exceptional divisor $E=\pi^{-1}(Z)$ of $\pi$.

{\rm (ii)} Furthermore, if $Z_0$ is closed in $Q$ then   $\tX$ contains a Zariski open set $X'$ isomorphic to 
$X$  (i.e. $X'$ is $k$-flexible) and with $X'\cap E$ being dense in $E$.
\ethm

\bproof The first statement follows from Theorem \ref{am.t1} in the case of $X=X_0$.
The second statement also follows from Theorem \ref{am.t1} since by the assumption the Serre Theorem A implies that
one can choose a regular function $f \in \C [Q]$ that vanishes
on $Z_0$ and equal to 1 on $D_0$ with any prescribed multiplicity.
\eproof

\brem\label{} (1) Consider statement (ii) of  Theorem \ref{kf.t1}.  Because of flexibility for any point $x \in X' \cap E$ and every nonzero vector $v\in T_xX'$
there exists a locally nilpotent vector $\sigma$ field on $X'$ for which $\sigma |_x =v$ \cite[Corollary 4.3]{AFKKZ}. In particular, this field
can be chosen transversal to $E$.

(2) Suppose that $Q_0=Q_1\times \C_{u}$ in the case of Theorem \ref{kf.t1} (i) and that the natural projection $Z_0 \to Q_1$ is \'etale.
Then the lift of the vector field $\partial / \partial u$ to $X_0'$ is also locally nilpotent and transversal to $X_0'\cap E$ at every point.
\erem

By the virtue of \cite[Proposition 4.14 and Remark 4.16]{AFKKZ} we have the following.

\bprop\label{kf.p1} Let $X$ be a $k$-flexible manifold and  let $\rho : X \to Q$ and $Q_0$  be as in Definition \ref{kf.d1}. Suppose that 
$x\in X$ is such that $\rho (x)\in Q_0, F=\rho^{-1} (\rho (x))$, and $V$ is a $k$-dimensional subspace of $T_xX$. 
Then there exists an automorphism $\alpha \in \Aut (X)$ such that $\alpha (x)=x$, $\alpha_* (T_xF)=V$,
and furthermore $\alpha_*$ transforms a given basis of $T_xF$ into a given basis of $V$.

\eprop

\section{Main theorems}

Now we are prepared for our main results.

\bthm\label{mt.t1} Let $X$ be a locally $k$-flexible manifold for $k\geq 2$, and $Z$ be a closed submanifold of $X$ of codimension at most $k$.
Suppose that $\pi : \tX \to X$ is the blowing up of $X$ along $Z$. Then $\tX$ is subelliptic. 

\ethm

\bproof  For the proof we can assume $X$ is $k$-flexible and by Remark \ref{kf.r1}(i) we can suppose that ${\rm codim}_X Z = k$.  Choose any point $z \in Z$ and any point $w \in \pi^{-1}(z)$. We need to construct a
family of sprays on $\tX$ of rank 1 that is dominating at $w$. By Proposition \ref{ff.p1} the phase flow of a complete vector field on $\tX$ can move $w$ in a general position
$\pi^{-1} (z)$, i.e. we can suppose that $w$ is general. 
Consider $\rho : X \to Q$ and $Q_0$ as in Definition \ref{kf.d1}, i.e. $X_0:=\rho^{-1}(Q_0)$ is 
 isomorphic
to $Q_0 \times \C^k$. Choose a morphism $\tau : X \to Q\times \C_u$ such that $\tau |_{X_0}= (\rho , \lambda )$ where
$\lambda : \C^k \to \C$ is a linear map.

By Theorem \ref{fm.t1}, applying an automorphism we can suppose that $\rho (z)\in Q_0$, $\tau|_Z : Z \to \tau (Z)=:Y$ is birational, $Y$ is a hypersurface smooth at $y=\tau (z)$,
and the projection $Y \to Q$ is smooth at $y$ (in particular, the vector field $\p /\p u$ is transversal to $Y$ at $y$. By Lemma \ref{am.l1} replacing this field $\p /\p u$
by an equivalent one $\delta$ we can extend it to $X$.
Let $X_0\simeq Q_0\times \C^k$ and $X_0'$ be
as in Theorem \ref{kf.t1}. Without loss of generality we can suppose that $Y\cap (Q_0\times \C)$ is smooth and is given in $Q_0\times \C$ by the zero locus of a
regular function $f$ on $Q_0\times \C$. Then $f\circ \pi$ yields a regular function on $X_0'$ whose zero locus may be viewed as a Zariski open subset $W$
of $ \pi^{-1}(Z)$ (and, moreover, this locus contains a Zariski open subset of $\pi^{-1}(z)$ for every $z \in Z$). 
Since $w$ is general we have $w \in W$.  Observe that because $X_0'\simeq X_0$ the field $\delta$ has a lift to a locally nilpotent vector 
field $\sigma$ on $X_0'$ which is transversal to $\pi^{-1}(Z)$ at $w$. By Proposition \ref{sas.p1} $\sigma$ extends to a spray of rank 1 on $\tX$
and the only thing we have to show that the vector $\sigma_w$ can be chosen general.

This follows from Proposition \ref{kf.p1} because we can transform the $\C^k$-fibration $\rho$ by some automorphism $\alpha$ into another $\C^k$-fibration
such that  $\alpha_* (\sigma_z)$ is a general vector. This yields the desired conclusion.

\eproof

As a corollary, we now give the proof of Theorem \ref{i.t1} in the introduction. 
\bproof
For the proof, we can assume that $X$ is stably flexible, i.e. $X$ is quasi-affine and $Y:=X\times \mathbb{C}^N$ is flexible for some positive integer $N$. Then, $Y$ is $N$-flexible, here we can choose $Q_0=X$.  Since the product of two flexible manifolds is again a flexible manifold, we can assume that $N\geq \dim (X)$. Consider $Z_1=Z\times \mathbb{C}^N$, then $Z_1$ is a smooth algebraic submanifold of $Y$, being of codimension $\leq \dim (X)\leq N$. Theorem \ref{mt.t1} implies that the blowup $\pi _1:\tY\rightarrow Y$ at $Z_1$ is algebraically subelliptic. Since $\tY=\tX\times \mathbb{C}^N$, it follows from the descent property for algebraic subellipticity that $\tX$ is itself algebraically subelliptic, as desired.     
\eproof

Here are some final remarks. 

\brem\label{mt.r2} (1) Note that $X=\SL_n$ is not contained in class $\mathcal{A}_0$ (or $\mathcal{A}$), i.e. it cannot be covered by open sets isomorphic
to $\C^N$ (where $N=\dim X$). Indeed,  $\SL_n$ is factorial since the ring of regular function on every simply connected algebraic group is a factorial domain (e.g, see \cite{Po}). 
Thus, if one assumes existence of an open subset $U \simeq \C^N$
such that $U \ne X$ then $D=X \setminus U$ must be a divisor because of affineness. Factoriality implies that $D=f^{-1}(0)$ for a regular function
on $X$. However this function must be constant on $U$ because of the fundamental theorem of algebra and, thus, on $X$.
A contradiction. 

(2)  
Let $H\subset \C_{u,v,\bar x}^{n+2}$ be a hypersurface
given by $uv=p(\bar x)$ in the case when the zero locus of $p$ is smooth connected. Then it is again factorial because of the Nagata lemma
(e.g., \cite{Ei}). Thus if $H$ is not isomorphic to $\C^{n+1}$ it does not belong to class $\cA$ by the same argument as before.



\erem

\end{document}